\documentclass{tac} 

\usepackage{hyperref}
\usepackage{amsmath,amssymb}
\usepackage[all]{xy}

\author{M. Menni}

\thanks{This work was supported by CONICET (Argentina), PIP 11220200100912CO and also by  the European Union’s Horizon 2020 research and innovation programme under the Marie Skłodowska-Curie grant agreement No 101007627.}

\address{Conicet and Universidad Nacional de La Plata, Argentina.}

\title[The least subtopos containing the discrete skeleton of $\Omega$]{The least subtopos containing the \\ discrete skeleton of $\Omega$}

\dedication{Dedicated to Pieter Hofstra.}

\copyrightyear{2023}


\eaddress{matias.menni@gmail.com}

\newcommand{\Set}{\mathbf{Set}}

\newcommand{\cat}[1]{\mathcal{#1}} 
\newcommand{\escat}[1]{\cat{#1}}

\newcommand{\Psh}[1]{\widehat{#1}}

\newcommand{\calC}{\ensuremath{\escat{C}}} 
\newcommand{\calE}{\ensuremath{\escat{E}}} 
\newcommand{\calF}{\ensuremath{\escat{F}}} 
\newcommand{\calL}{\ensuremath{\escat{L}}} 
\newcommand{\calS}{\ensuremath{\escat{S}}} 

\newcommand{\opCat}[1]{\ensuremath{{#1}^{\mathrm{op}}}}

\newcommand{\OmegaBase}{\mathit{2}}

\begin{document}

\maketitle

\begin{abstract}
Let ${p : \calE\rightarrow \calS}$ be  a pre-cohesive geometric morphism.
We show that the least subtopos of $\calE$ containing both the subcategories ${p^* : \calS \rightarrow \calE}$ and ${p^! : \calS \rightarrow \calE}$ exists, and that it coincides with the least subtopos containing ${p^* \OmegaBase}$, where  ${\OmegaBase}$ denotes the subobject classifier of $\calS$.
\end{abstract}


\section{Introduction}

At least since his PhD  \cite{Hofstra2003}, and throughout his career until \cite{FunkHofstra2022}, Topos Theory was one of Pieter's many interests.
Each time we met we spoke about it. Certainly a lot when Pino Rosolini invited him to Genova while I was there in 2002, and then also a little each time we met at CT meetings.  He was a very witty conversationalist so it was also a pleasure to discuss other things. Still,  I would have liked to tell him about the result in the present paper because it concerns toposes `of spaces' which, as explained in \cite[Theorem~1]{Lawvere07}, have a canonical Homotopy Theory, another one of Pieter's interests \cite{HofstraWarren2013}.

More specifically, this paper is about the Dimension Theory outlined in \cite{Lawvere91}.  
We recall some of the basic definitions and refer to  \cite{MMlevelEpsilon, Menni2019a} for  additional information.
A {\em level} of a topos $\calE$  is just an essential subtopos of $\calE$.
If ${l : \calL \rightarrow \calE}$ is a level of $\calE$, then the left adjoint to the inverse image functor  ${l^* : \calE \rightarrow \calL}$ is denoted by ${l_!}$. Following \cite{Lawvere91} we picture a level of $\calE$ as a `dimension' and ${l_! : \calL \rightarrow \calE}$ as the full subcategory  consisting of the objects $X$ of $\calE$ such that `${\dim X \leq l}$'. For any $X$ in $\calE$, the counit ${l_! (l^* X) \rightarrow X}$ is called the { \em ($l$)-skeleton} of $X$. 
Abusing the terminology a bit, we sometimes call the object ${l_! (l^* X)}$ the { \em ($l$)-skeleton} of $X$. 
In other words, we sometimes confuse a skeleton with its domain.

Levels  may be partially ordered as subtoposes so, in this way, each topos determines its poset `of dimensions'.
Also, for levels $l$, $m$ of a given topos, we say that $m$ is {\em way above} $l$ if $m$ is above $l$ as subtoposes and, moreover,  ${l_!}$ factors through $m_*$. The {\em Aufhebung} of a level $l$ is the least level that is way above $l$. 

\begin{remark}\label{RemDefAufhebung} 
Notice that it also makes sense to say that a subtopos is way above some  level $l$. Hence, we may consider the problem of finding the least subtopos that is way above  $l$. Of course, if this least subtopos is essential then it is the Aufhebung of $l$. The resulting `weak Aufhebung' problem seems interesting because, on the one hand, little is known about  conditions guaranteeing the existence of the Aufhebung of an arbitrary level of an elementary topos and, on the other hand, for several important toposes `of combinatorial spaces', such as simplicial sets, every subtopos is essential.
\end{remark}

Remark~\ref{RemDefAufhebung} and other well-known results in Topos Theory suggest the following.

\begin{definition}\label{DefEnvelope}{\em 
For a topos   $\calE$  and a functor ${F : \calC \rightarrow \calE}$, the {\em envelope} of $F$ is, if it exists, the  least subtopos ${f : \calF \rightarrow \calE}$ such that $F$ factors through the direct image functor ${f_* : \calF \rightarrow \calE}$.}
\end{definition}

Roughly speaking, the envelope of $F$ is the  least subtopos of $\calE$ through which $F$ factors.
As suggested above, given a level ${l : \calL \rightarrow \calE}$, the `weak Aufhebung' of $l$ is the envelope of  the join of the subcategories ${l_! : \calL \rightarrow \calE}$ and ${l_* : \calL\rightarrow \calE}$. 
The main result of the paper shows that it exists in certain cases. 
Before stating the result precisely it is relevant to mention other examples of envelopes.

It makes perfect sense to consider the envelope of an object in a topos.
Indeed, we may define it as the envelope of the obvious functor from the terminal category to the  topos.
The envelope of any object in a topos exists by  \cite[A4.5.15]{elephant}.

\begin{example}[The least dense subtopos.]\label{ExLeastDenseSubtopos}
The envelope of the initial object coincides with the subtopos for the double-negation topology \cite[paragraph after A4.5.20]{elephant}.
\end{example}

\begin{example}[Weak generation.]
A topos is said to be {\em weakly generated} by an object if the envelope of that  object is the whole topos. For example,  \cite[A2.3.9]{elephant} shows that every topos is weakly generated by its subobject classifier; but toposes may be weakly generated by smaller objects. The topological topos is weakly generated by 2 \cite[Example~{1.6}]{Menni2019}. Sufficiently Cohesive toposes over Boolean bases are weakly generated by objects with a unique point \cite[Corollary~{6.3}]{Menni2019}.
\end{example}

If ${p : \calE\rightarrow \calS}$ is a geometric morphism and  $\OmegaBase$ is the subobject classifier of $\calS$, then we may consider the envelope of  ${p^* \OmegaBase}$.

\begin{example}[The least pure subtopos.]
If ${p : \calE \rightarrow \calS}$ is locally connected and bounded, then \cite[Proposition~{9.2.10}]{BungeFunkBook}  shows that the envelope of  ${p^* \OmegaBase}$ is the smallest pure subtopos of $\calE$ and it is also locally connected over $\calS$.
\end{example}

For a level ${l : \calL \rightarrow \calE}$  we may consider the envelope of  $l_! \OmegaBase$ where $\OmegaBase$ is the subobject classifier of $\calL$, and we may also consider the envelope of  the (domain of the) $l$-skeleton of the subobject classifier of $\calE$. Occasionally, these two subtoposes coincide; in particular, this is the case if the inverse image functor $l^*$ preserves the subobject classifier \cite[A4.5.8]{elephant}.

\begin{example}[The `weak Aufhebung' of level $-\infty$.]\label{ExWeakAofMinusInfty}
Level $-\infty$ is the essential subtopos  ${ \mathbf{0} \rightarrow \calE}$ where $\mathbf{0}$ is the topos with exactly one object. 
The inverse image functor trivially preserves the subobject classifier and the skeleton of any object in $\calE$ is initial. So, by Example~\ref{ExLeastDenseSubtopos},  the envelope of  the ${-\infty}$-skeleton of the subobject classifier coincides with the subtopos of sheaves for the double-negation topology. 
\end{example}

Let ${p : \calE \rightarrow \calS}$ be a geometric morphism. The map $p$ is {\em hyperconnected} if  both the unit and counit of ${p^* \dashv p_*}$ are monic. The map $p$ is {\em local} if its direct image functor has a fully faithful right adjoint, usually denoted by ${p^!}$. In this case, the string ${p^* \dashv p_* \dashv p^! : \calS \rightarrow \calE}$ is a level of $\calE$ called the {\em centre} of $p$, and, for any $X$ in $\calE$, the counit ${p^* (p_* X ) \rightarrow X}$ is the skeleton of $X$ (determined by the centre of $p$).
 Also, $p$ is {\em pre-cohesive} if it is local, hyperconnected and the inverse image functor has a finite-product preserving left adjoint, usually denoted by $p_!$  \cite{LawvereMenni2015}. In this case, we picture the objects of $\calE$ as `spaces' and ${p^* : \calS \rightarrow \calE}$ as the full subcategory of `discrete spaces'. Also,  ${p_!}$  is thought of as sending a space $X$ to the discrete space ${p_! X}$ of connected components. Then, for every space $X$, the (monic) skeleton ${p^* (p_* X) \rightarrow X}$ is called the {\em discrete skeleton} of $X$.
Allowing ourselves to identifying `discrete' with `0-dimensional', it makes sense to say that  the centre of the pre-cohesive $p$ is {\em level 0}.

We prove that, for a pre-cohesive ${p : \calE \rightarrow \calS}$, the `weak Aufhebung' of level~0 exists and, moreover, it coincides with the least subtopos containing the 0-skeleton of $\Omega$. Notice that this is analogous to  the fact that, if $\calS$ is Boolean, then level~0 is the Aufhebung of level~{$-\infty$}, as already observed in  \cite[Corollary~{4.5}]{LawvereMenni2015} and Example~\ref{ExWeakAofMinusInfty}.

We stress that for an arbitrary geometric morphism ${p : \calE \rightarrow \calS}$ towards a topos $\calS$ with subobject classifier $\mathit{2}$, the envelope of ${p^* \mathit{2}}$ need not be a level of $\calE$, even if $p$ is locally connected.

\begin{example}[The smallest pure subtopos of ${\mathrm{Sh}(\mathbb{R}^2)}$  is  not essential.]
The present example is due to Jon Funk who patiently explained it to me in a private communication.
Let ${U \subseteq \mathbb{R}^2}$ be a pure open subet and ${r \in U}$.
Then ${U - \{r\} \subseteq \mathbb{R}^2}$ is also pure and open (and strictly smaller than $U$, of course).
Hence, the smallest pure subtopos of ${\mathrm{Sh}(\mathbb{R}^2)}$ is not open.
Also,  \cite{Hemelaer2022} implies that the levels of a topos of sheaves on a Hausdorff space are precisely its open subtoposes.
Therefore, the smallest pure subtopos of ${\mathrm{Sh}(\mathbb{R}^2)}$  cannot be essential, for otherwise it would be open.
\end{example}

In contrast, if $\calE$ is a topos such that every subtopos is a level, then, trivially, every envelope is a level. In such toposes there is an endofunction on the  poset of levels of $\calE$ that sends a level $l$ to the envelope of  the $l$-skeleton of the subobject classifier of $\calE$. As a source of examples consider the topos of presheaves on any small category with split-epi/split-mono factorizations and such that every object has a finite set of subobjects \cite[Proposition~{2.5}]{Menni2019a}. In particular, we may consider the topos of simplicial sets, or the classifier of non-trivial Boolean algebras. We don't know if any of these functions coincides with the Aufhebung of the respective topos.

\section{Internal and stable orthogonality}

Let $\calE$ be a topos with subobject classifier $\Omega$ and let ${f : W \rightarrow X}$ be a morphism in $\calE$.
We say that $f$ is {\em orthogonal} to an object $Z$ if for every ${g : W \rightarrow Z}$ there exists a unique ${g' : X \rightarrow Z}$ such that ${g' f = g}$. Alternatively, $f$ is orthogonal to $Z$ if and only if ${\calE(f, X) : \calE(X , Z)  \rightarrow \calE(W, Z)}$ is an isomorphism. 
We say that $f$ is {\em internally orthogonal} to $Z$ if ${Z^f : Z^X \rightarrow Z^W}$ is an isomorphism.

\begin{remark}\label{Remark} It is well-known that ${\Omega^{(\_)} : \opCat{\calE} \rightarrow \calE}$  reflects isomorphisms. In other words, a map $f$ is internally orthogonal to $\Omega$ if and only if $f$ is an isomorphism.
\end{remark}

The map $f$ is {\em stably orthogonal} to $Z$ if for every pullback as below
$$\xymatrix{
U \ar[d]_-{f'} \ar[r] & W \ar[d]^-f \\
V \ar[r] & X
}$$
the map $f'$ is orthogonal to $Z$.

\begin{lemma}\label{LemStablyImpliesInternallyOrtho} If $f$ is stably orthogonal to $Z$, then it is internally so.
\end{lemma}
\begin{proof}
Simply observe that  the map $f$ is internally orthogonal to $Z$ if and only if,  for every object $Y$,
${f \times Y : W \times Y \rightarrow X \times Y}$ is orthogonal to $Z$ .
\end{proof}

The proof that the least subtopos containing an object exists  \cite[Proposition~{A4.5.15}]{elephant}
 involves  the following.

\begin{lemma}\label{LemDenseMonosWRTobject} 
If ${\calE_j \rightarrow \calE}$ is  the least subtopos containing the object ${Z}$ in $\calE$ then, 
a subobject ${u : U \rightarrow X}$ in $\calE$ is $j$-dense   if and only if $u$ is stably orthogonal to $Z$.
\end{lemma}

Combining Lemmas~\ref{LemStablyImpliesInternallyOrtho} and \ref{LemDenseMonosWRTobject} we obtain the next auxiliary fact.

\begin{lemma}\label{LemDenseForEnvelopeOfObjectImpliesIO} If ${\calE_j \rightarrow \calE}$ is the least subtopos containing the object ${Z}$ in $\calE$, then every $j$-dense subobject is internally orthogonal to $Z$. 
\end{lemma}

Let $\calS$ be another topos with subobject classifier denoted by $\mathit{2}$.

\begin{lemma}\label{LemOrthoTo2Abstract}  If  ${R : \calS \rightarrow \calE}$ is a full and faithful functor with a finite-product preserving left adjoint  ${L : \calE \rightarrow \calS}$, then   the following are equivalent:
\begin{enumerate}
\item The map $f$ is internally orthogonal to ${R \mathit{2}}$.
\item The map ${L f: L W \rightarrow L X}$ is an isomorphism  in $\calS$.
\item For every $A$ in $\calS$, $f$ is internally orthogonal to ${R A}$.
\end{enumerate}
\end{lemma}
\begin{proof}
As $L$ preserves finite products by hypothesis, the canonical  ${ R(A^{L Y})  \rightarrow (R A)^Y}$ is an isomorphism natural in both $A$ and $Y$.
In particular we have
$$\xymatrix{
(R \mathit{2})^X \ar[d]_-{\cong} \ar[rr]^-{(R \mathit{2})^f} && (R \mathit{2})^W \ar[d]^-{\cong} \\
R( \mathit{2}^{L X}) \ar[rr]_-{R(\mathit{2}^{L f})} && R (\mathit{2}^{L W})
}$$
so, if the top map is an isomorphism, then the bottom map is an isomorphism.
As $R$ is fully faithful, ${\OmegaBase^{L f}}$ is an isomorphism. That is, ${L f}$ is internally orthogonal to $\OmegaBase$.
So ${L f}$ is an isomorphism by Remark~\ref{Remark}.

To prove that the second item implies the third consider the diagram
$$\xymatrix{
(R A)^X \ar[d]_-{\cong} \ar[rr]^-{(R A)^f} && (R A)^W \ar[d]^-{\cong} \\
R( A^{L X}) \ar[rr]_-{R(A^{L f})} && R (A^{L W})
}$$
where the vertical maps are iso by adjointness and the bottom map is iso by hypothesis.
It follows that the top map is also an iso, showing that $f$ is internally orthogonal to ${R A}$.

Finally, the third item trivially implies the first.
\end{proof}

\section{The `weak Aufhebung' of level~0  of a pre-cohesive topos}

Let ${p : \calE \rightarrow \calS}$ be a pre-cohesive geometric morphism.
The centre of the local $p$ is called {\em Level~0}.
 Let ${\OmegaBase}$ be the subobject classifier  of $\calS$ but we stress  that we are not assuming that $\calS$ is Boolean.

\begin{lemma}\label{LemCohesiveOrthoTo2} For every ${f : W \rightarrow X}$ in $\calE$, the following are equivalent:
\begin{enumerate}
\item The map $f$ is internally orthogonal to ${p^* \mathit{2}}$.
\item The map ${p_! f : p_! W \rightarrow p_! X}$ is an isomorphism in $\calS$.
\item For every $A$ in $\calS$, $f$ is internally orthogonal to ${p^* A}$.
\end{enumerate}
\end{lemma}
\begin{proof}
Follows from Lemma~\ref{LemOrthoTo2Abstract} applied to the adjunction ${p_! \dashv p^*}$.
\end{proof}

In other words, $f$ is internally orthogonal to discrete spaces if and only if $f$ induces an isomorphism at the level of connected components.
Notice that the equivalence between the first and third items of Lemma~\ref{LemCohesiveOrthoTo2} appears to be a variant of \cite[Proposition~{9.2.14}]{BungeFunkBook}.
We stress  that  monomorphisms satisfying the equivalent conditions of Lemma~\ref{LemCohesiveOrthoTo2} need not be stable under pullback.

\begin{example} 
Consider the topos ${\Psh{\Delta_1}}$ of reflexive graphs and the pre-cohesive geometric morphism ${p : \Psh{\Delta_1} \rightarrow \Set}$. Let  ${X = \cdot \rightrightarrows \cdot}$ be the  reflexive graph with two nodes and two (non-identity) parallel edges. Pick one of the non-identity edges and consider the associated subobject ${u}$ of ${X}$. Its intersection with the subobject determined by the other edge is the discrete set of nodes of $X$. So $u$ is a monomorphism such that ${p_! u}$ is an isomorphism, but it is not stably orthogonal to ${p^* 2}$. Notice also that in this example, ${p_* u}$ is  an isomorphism. 
\end{example}

\begin{lemma}\label{LemCohesiveStablyOrthoTo2} Let ${\calE_j \rightarrow \calE}$ be least subtopos containing   ${p^* \mathit{2}}$. A monic ${u : U \rightarrow X}$ in $\calE$ is $j$-dense if and only if, for every pullback square 
$$\xymatrix{
V  \ar[d]_-v  \ar[r] & U \ar[d]^-u \\
Y \ar[r] & X
}$$
in $\calE$,  ${p_! v}$ is an isomorphism in $\calS$. Hence, if $u$ is $j$-dense then ${p_* u}$ is an isomorphism.
\end{lemma}
\begin{proof}
Lemma~\ref{LemCohesiveOrthoTo2} implies that the pullback condition in the statement holds if and only if $u$ is stably orthogonal to ${p^* \mathit{2}}$.
So the first part of the result follows from Lemma~\ref{LemDenseMonosWRTobject}. To prove the second part of the statement observe that the square below 
$$\xymatrix{
p^* (p_* U) \ar[d]_-{p^* (p_* u)} \ar[r]^-{\beta} & U \ar[d]^-u \\
p^* (p_* X) \ar[r]_-{\beta} & X
}$$
is a pullback because $p$ is hyperconnected so ${p_! (p^* (p_* u)) : p_! (p^* (p_* U)) \rightarrow p_!(p^* (p_* X))}$ is an isomorphism in $\calS$.
That is, ${p_* u}$ is an isomorphism.
\end{proof}

\begin{example}
Let  ${p : \Psh{\calC} \rightarrow \Set}$ be the canonical geometric morphism from the topos of presheaves on the small category $\calC$.
The map $p$ is not necessarily pre-cohesive but it is locally connected so $p^*$ has a left adjoint $p_!$ and \cite[Example~{9.2.12}]{BungeFunkBook} implies that a sieve $R$ on an object $C$ of $\calC$ is covering for the least subtopos containing ${p^* \OmegaBase}$  if and only if ${f^* R}$ is connected for every map $f$ with codomain $C$ in $\calC$. This  is  essentially the pullback condition of Lemma~\ref{LemCohesiveStablyOrthoTo2} but restricted to the site. So this example illustrates that lemma and suggests a possible generalization to locally connected, but not necessarily pre-cohesive, maps.
\end{example}

\begin{theorem}\label{ThmWAofL0} Let $\calS$ be a topos with subobject classifier $\OmegaBase$. If ${p : \calE \rightarrow \calS}$ is pre-cohesive, then  the  least subtopos of $\calE$ containing the subcategories ${p^* , p^! : \calS \rightarrow \calE}$ exists and it coincides with the least subtopos of $\calE$ containing ${p^* \OmegaBase}$.
\end{theorem}
\begin{proof}
 Let ${\calE_j \rightarrow \calE}$ be  the least subtopos containing  ${p^* \OmegaBase}$ and let $A$ be an object in $\calS$.
 Lemma~\ref{LemCohesiveStablyOrthoTo2} easily implies that every $j$-dense mono is orthogonal to ${p^! A}$. 
 That is, ${p^! A}$ is a $j$-sheaf and therefore, ${\calE_j \rightarrow \calE}$ is above the centre of $p$.

It remains to prove that the inclusion ${p^* : \calS \rightarrow \calE}$ factors through the inclusion ${\calE_j \rightarrow \calE}$. 
By Lemma~\ref{LemDenseForEnvelopeOfObjectImpliesIO}, $j$-dense subobjects are internally orthogonal to ${p^* \mathit{2}}$  and therefore orthogonal to ${p^* A}$ by 
Lemma~\ref{LemCohesiveOrthoTo2}. In other words, ${p^* A}$ is a $j$-sheaf.

We have proved that the subtopos ${\calE_j \rightarrow \calE}$ is way-above level~0.
It must be the least one because any subtopos of $\calE$ that is way-above level~0 must contain ${p^* \mathit{2}}$.
\end{proof}

\end{document}